\documentstyle[11pt]{article}
\newlength{\spacing}
\newcommand{\doublespace}{\setlength{\baselineskip}{1.5\spacing}}

\newcommand{\eq}[1]{\begin{equation} #1 \end{equation}}

\newtheorem{thm}{Theorem}[section]

\newtheorem{lem}[thm]{Lemma}
\newtheorem{prop}[thm]{Proposition}

\def\sec{\setcounter{equation}{0}}

\def\lam{\lambda }

\def\ka{\kappa}
\def\lab{\label }

\def\rar{\to}

\def\inft{\infty}

\def\ep{\epsilon}
\def\sg{\sigma}

\def\cP{{\cal P}}

\def\cR{{\Re}}

\def\cX{{\cal X}}

\def\today{\ifcase\month\or
  January\or February\or March\or April\or May\or June\or
  July\or August\or September\or October\or November\or December\fi
  \space\number\day, \number\year}

\begin{document}
\begin{titlepage}
\begin{center}
{\bf A Strong Law for the Largest Nearest-Neighbor Link on
Normally Distributed Points} \\
\vspace{0.20in} by \\
\vspace{0.2in} Bhupender Gupta \\
Department of Mathematics, Indian Institute of Technology,
Kanpur 208016, India \\
\vspace{0.1in} and \\
\vspace{0.2in} {Srikanth K. Iyer \footnote{Corresponding Author.
email:skiyer@math.iisc.ernet.in}}\\
Department of Mathematics,
Indian Institute of Science, Bangalore 560012, India. \\
\end{center}
\vspace{0.2in}
%
\sloppy
\begin{center} {\bf Abstract} \end{center}

\begin{center} \parbox{4.8in}
{Let $n$ points be placed independently in $d-$dimensional space
according to the standard $d-$dimensional normal distribution.
Let $d_n$ be the longest edge length for the nearest neighbor
graph on these points. We show that
\[\lim_{n \rar \infty} \frac{\sqrt{\log n} \; d_n}{\log \;
\log n} =
\frac{d}{\sqrt{2}}, \qquad d \geq 2, \mbox{ a.s.} \] } \\
\vspace{0.4in} \today
\end{center}

\vspace{0.5in}
{\sl AMS 1991 subject classifications}: \\
\hspace*{0.5in} Primary:   60JD05, 60G70\\
\hspace*{0.5in} Secondary:  05C05, 90C27

\end{titlepage}
\doublespace
\section{Introduction and main results\lab{s1}}
\sec

In this paper we prove a strong law result for the largest
nearest neighbor distance of points distributed according to a
standard normal distribution in $\cR^d.$ Throughout this paper we
will assume that $d \geq 2.$

Let $X_1,X_2,\ldots$ be independent and identically distributed
random variables distributed according to the standard
multivariate normal distribution in $\cR^d$. Let $\phi(x),$ $x
\in \cR^d$, denote the standard multivariate normal density,
\[ \phi(x) = (2\pi)^{-d/2}\exp(-\|x\|^2/2), \]
where $\| \cdot \|$ is the Euclidean $(\ell_2)$ norm on $\cR^d.$
Let $R = \| X \|.$ Then the probability density function of $R$ is
given by
\begin{equation}
f_{R}(r) = A_d e^{- r^2/2}r^{d-1},\qquad 0<r<\inft, d\geq 2
\lab{e0}
\end{equation}
where $A_d = \frac{(2\pi)^{-d/2}}{d}.$

The basic object of study will be the graphs $G_n$ with vertex
set $\cX_n = \{X_1,X_2,\ldots ,X_n\}$, $n=1,2,\ldots.$ Edges of
$G_n$ are formed by connecting each of the vertices in $\cX_n$ to
its nearest neighbor. The longest edge of the graph $G_n$ is
denoted by $d_n$. We shall refer to $G_n$ as the nearest neighbor
graph (NNG) on $\cX_n$ and to $d_n$ as the largest nearest
neighbor distance (LNND).

The largest nearest neighbor link has been studied in the context
of computational geometry (see Dette and Henze (1989) and Steele
and Tierney (1986)) and has applications in statistics, computer
science, biology and the physical sciences. Appel and Russo
(1997) proved strong law results for a graph on uniform points in
the $d-$ dimensional unit cube. Penrose (1999) extended this to
general densities having compact support. Penrose (1998) proved a
weak law result for $d_n$ on normally distributed points, which
states that $\sqrt{(2\log n)} d_n - b_n$ converges weakly to the
Gumbel distribution, where $b_n \sim (d-1) \log \log n.$ $a_n
\sim b_n$ implies that $a_n/b_n$ converges to a constant as $n
\rar \infty.$ In what follows we will write $\log_2 n$ for $\log
\log n$. The above result is also shown to be true for the
longest edge of the minimal spanning tree. We are not aware of
strong law results for the LNND for graphs whose vertices are
distributed according to densities with unbounded supports for $d
\geq 2.$ For a detailed description of Random Geometric Graphs,
their properties and applications, we refer the reader to Penrose
(2003) and references therin.

It is often easier to study the graph $G_n$ via the NNG $P_n$ on
the set $\cP_n = \{X_1,X_2,\ldots ,X_{N_n}\}$, where $N_n$ are
Poisson random variables with mean $n.$ Then $\cP_n$ is an
inhomogeneous Poisson point process with intensity function
$n\phi(\cdot).$ Note that the graphs $G_n$ and $P_n$ are coupled,
since the first $\min(n,N_n)$ vertices of the two graphs are
identical. We also assume that the random variables $N_n$ are
non-decreasing, so that $P_1 \subset P_2 \subset P_3 \cdots.$ We
now state our main result.

\begin{thm} Let $d_n$ be the LNND of the NNG $\; G_n,$ defined on
the collection $\cX_n$ of $n$ points distributed independently
and identically according to the standard normal distribution in
$R^d,$ $d \geq 2.$ Then,
\begin{equation} \lim_{n \rar \infty } \frac{\sqrt{\log n}\;
d_n}{\log_2 n} = \frac{d}{\sqrt{2}}, \qquad d \geq 2, \mbox{ a.s.}
\label{e1} \end{equation} \label{t1} \end{thm}

\vskip0.5truein
\section{Proofs and supporting results}
\sec

For any $x \in \cR^d,$ let $B(x,r)$ denote the open ball of
radius $r$ centered at $x.$ Let
\eq{ I(x,r) = \int_{B(x,r)} \phi(y) \; dy. \lab{e2} }
For $\rho > 0,$ define $I(\rho,r) = I(\rho e, r),$ where $e$ is
the $d-$dimensional unit vector $(1,0,0,\ldots, 0).$ Due to the
radial symmetry of $\phi(x),$ $I(x,r)=I(\| x \| ,r).$ The
following Lemma (see Penrose (2003), Lemma 8.8) provides large
$\rho$ asymptotics for $I(\rho,r).$

\begin{lem} Let $(\rho_n)_{n \geq 1}$ and $(r_n)_{n \geq 1}$ be
sequences of positive numbers such that $\rho_n \rar \infty,$
$r_n \rar 0,$ and $r_n \rho_n \rar \infty,$ as $n \rar \infty.$
Then,
\eq{ I(\rho_n,r_n) \sim (2 \pi)^{-1/2} r_n^d \exp(\rho_nr_n -
\rho_n^2) (\rho_n r_n)^{-(d+1)/2}, \qquad \mbox{ as } n \rar
\infty. \lab{e3}} \lab{l1} \end{lem}
%
%
%
In order to prove strong law results for the LNND for graphs with
densities having compact support, one covers the support of the
density using an appropriate collection of concentric balls and
then shows summability of certain events involving the
distribution of the points of $\cX_n$ on these balls. The results
then follow by an application of the Borel-Cantelli Lemma. The
asymptotic behaviour of the LNND depends on the (reciprocal of
the) infimum of the density, since it is in the vicinity of this
infimum, points will be sparse and hence be farthest from each
other (see for example Penrose (1999)). In case of densities
having unbounded support, the region to be covered must be
determined first. The following Lemma gives us the regions of
interest when the points are normally distributed.

For any $c \in \cR$ fixed, and large enough $n$, define
\eq{ R_n(c) = \sqrt{ 2 \log n + (c+d-2) \log_2 n + 2 \log(A_d)},
\lab{e4a}}
where $A_d$ is as in (\ref{e0}). Let $A^c$ denote the complement
of set $A.$ Let $U_n(c)$ be the event $\cX_n \subset B(0,R_n(c))$
and $V_n(c)$ denote the event that at least one point of $\cX_n$
lies in $B^c(0,R_n(c)).$ $a_n \stackrel{>}{\sim} b_n$ implies
that $a_n > c_n$ for some sequence $c_n$ and $c_n \sim b_n.$

\begin{lem} For $c > 2$, $P[U_n^c(c) \mbox{ i.o. }] = 0,$ and
for $c < 0$, $P[V_n^c(c) \mbox{ i.o. }] = 0$ a.s. The result is
also true with $\cX_n$ replaced by $\cP_{\lam_n}$ provided $\lam_n
\sim n.$ \lab{l2}
\end{lem}

Thus for almost all realizations of the sequence $\{\cX_n\}_{n
\geq 1}$, all points of $\cX_n$ will lie within the ball
$B(0,R_n(c))$ $c > 2,$ eventually and for $c < 0$, there will be
at least one point of $\cX_n$ in $B^c(0,R_n(c))$ eventually.

{\bf Proof of Lemma~\ref{l2}. } As $R_n \rar \infty,$ note that
\eq{ 1 - I(0,R_n) := \int_{R_n}^{\infty} A_d e^{-R_n^2/2}
R_n^{d-1} \sim A_d R_n^{d-2}e^{-R_n^2/2} \lab{e0l2}}
Hence,
\eq{ P[U_n^c(c)] = 1 - (I(0,R_n(c)))^n \sim A_d n
R_n(c)^{d-2}\exp(-R_n(c)^2/2). \lab{e3l2}}
Let $n_k$ be the subsequence $a^k,$ with $a > 1,$ and consider
\begin{eqnarray}
P[\cup_{n= n_k}^{n_{k+1}}U_n^c(c)] & \leq &
P[\mbox{at least one vertex of ${\cal{X}}_{n_{k+1}}$ is in $B^c(0,R_{n_k}(c))$}] \nonumber\\
& = & 1- (I(0,R_{n_k}(c)))^{n_{k+1}} \nonumber \\
& \sim & A_d n_{k+1} R_{n_k}^{d-2}(c)e^{-R_{n_k}^2(c)/2} \sim
k^{-c/2}. \lab{e1l2}
\end{eqnarray}
Thus the above probability is summable for $c>2$  and the first
part of Lemma~\ref{l2} follows from the Borel-Cantelli Lemma.
Again, using (\ref{e0l2}) and the inequality $1-x \leq \exp(-x),$
we get
\[ P[V_n^c(c)] = (I(0,R_n(c))^n \leq \exp(-A_d n R_n^{d-2}
e^{-R_n^2/2}). \]
Let $n_k$ be as above.
\begin{eqnarray}
P[\cup_{n= n_k}^{n_{k+1}}V_n^c(c)] & \leq &
P[\cX_{n_k} \subset B(0,R_{n_{k+1}}(c))] \nonumber \\
& \leq & \exp(-A_d n_k R_{n_{k+1}}^{d-2} e^{-R_{n_{k+1}}^2/2})
\nonumber \\
& \sim & \exp(- \mbox{ constant }k^{-c/2}),
\end{eqnarray}
which is summable for all $c < 0.$ This proves the second part of
Lemma~\ref{l2}.

If $\cX_n$ is replaced by $\cP_{\lam_n}$, then
\begin{eqnarray*}
P[U_n^c(c)] & = & \exp(-\lam_n (1-I(0,R_n(c)))) \\
& \stackrel{<}{\sim} & \lam_n (1 - A_d
R_n(c)^{d-2}\exp(-R_n(c)^2/2)) \sim n A_d
R_n(c)^{d-2}\exp(-R_n(c)^2/2),
\end{eqnarray*}
which is same as the right hand side of (\ref{e3l2}). Similarly,
one can show that $P[V_n(c)]$ has the same asymptotic behavior as
in the case of $\cX_n.$ Thus the results stated for $\cX_n$ also
hold for $\cP_{\lam_n}$.

\begin{prop} Let $t > \frac{d}{\sqrt{2}},$ and let
$r_n(t) = \frac{t \log_2 n}{\sqrt{\log n}}.$ Then with probability
$1$, $d_n < r_n(t)$ for all large enough $n.$ \lab{prop1}
\end{prop}
{\bf Proof. } Pick $u,t$ such that $(2d + c -2)/2\sqrt{2} < u <
t,$ and $\ep > 0$ satisfying
\[ \ep + u < t. \]
Let $c>2.$ From Lemma~\ref{l2}, $\cX_n \subset B(0,R_n(c))$ a.s.
for all large enough $n.$ For $m=1,2, \ldots,$ let $\nu(m) = a^m,$
for some $a > 1.$ Let $\ka_m$, (the covering number) be the
minimum number of balls of radius $r_{\nu(m)}(\ep)$ required to
cover the ball $B(0,R_{\nu(m+1)}(c)).$ From Lemma 2.1, Penrose
(1999), we have
\eq{ \ka_m \stackrel{<}{\sim} \left(\frac{m}{\log \; m}\right)^d.
\lab{p4}}
Consider the deterministic set $\{x_1^m,\ldots ,x_{\ka_m}^m\} \in
B(0,R_{\nu(m+1)}(c)),$ such that $B(0,R_{\nu(m+1)}(c)) \subset
\cup_{i=1}^{\ka_m} B(x_i^m,r_{\nu(m)}(\ep)).$

Given $x \in \Re^d,$ let $A_m(x)$ denote the annulus
$B(x,r_{\nu(m+1)}(u))\setminus B(x,r_{\nu(m)}(\ep)),$ and let
$F_m(x)$ be the event such that no vertex of ${\cal{X}}_{\nu(m)}$
lies in $A_m(x),$ i.e.
\begin{equation}
F_m(x) = \{{\cal{X}}_{\nu(m)}[A_m(x)]=0 \}\label{fm}
\end{equation}
Since,
\begin{eqnarray}
P[X_1 \in A_m(x)] & = & \int_{A_m(x)} \phi(y) \; dy \nonumber \\
& \geq & \int_{A_m(R_{\nu(m+1)}(c))} \phi(y) \; dy \nonumber \\
& = & I(R_{\nu(m+1)}(c), r_{\nu(m)}(u)) - I(R_{\nu(m+1)}(c),
r_{\nu(m)}(\ep)) \label{annulus pro}
\end{eqnarray}
from Lemma~\ref{l1}, we get
\begin{eqnarray*}
\lefteqn{P[X_1 \in A_n(x)]  \stackrel{>}{\sim} C_d
e^{-R_{\nu(m+1)}^2(c)/2}(R_{\nu(m+1)}(c))^{-\frac{d+1}{2}} \cdot} \\
& & \cdot \left(e^{R_{\nu(m+1)}(c)
r_{\nu(m)}(u)}(r_{\nu(m)}(u))^{\frac{d-1}{2}} - e^{R_{\nu(m+1)}(c)
r_{\nu(m)}(\ep)}(r_{\nu(m)}(\ep))^{\frac{d-1}{2}}\right) := q_m.
\end{eqnarray*}
Substituting the values of $R_{\nu(m+1)}(c)$ and
$r_{\nu(m)}(\cdot)$ in $q_m$, we get
\begin{eqnarray}
q_m & \sim & \frac{(\log \; m)^{(d-1)/2}}{a^{m+1}
m^{d+\frac{c}{2} - 1 - u \sqrt{2}}} \label{qn}
\end{eqnarray}
Hence,
\eq{ P[F_m(x)] \stackrel{<}{\sim}  (1-q_m)^{\nu(m)} \leq
\exp(-\nu(m) q_m ) \sim \exp\left\{ - C \frac{(\log \;
m)^{(d-1)/2}}{m^{d+\frac{c}{2} - 1 - u \sqrt{2}}} \right\},
\lab{p5}}
where $C$ is some constant. Set $G_m =
\cup^{\ka_m}_{i=1}F_m(x_i^m).$
\begin{eqnarray}
P[G_m] & \stackrel{<}{\sim} & \left(\frac{m}{\log \; m}\right)^d
\exp\left\{ - C \frac{(\log \; m)^{(d-1)/2}}{m^{d+\frac{c}{2} - 1
- u \sqrt{2}}} \right\},
\end{eqnarray}
which is summable in $m$ for all $u > (2d + c - 2)/2\sqrt{2}.$ By
Borel-Cantelli, $G_m$ occurs only for finitely many $m$ a.s.

Pick $n$,and take $m$ such that $a^m \leq n \leq a^{m+1}.$ If
$d_n \geq r_n(t)$, then there exists an $X \in \cX_n$ such that
$\cX_n(B(X,r_n(t))\setminus \{X\}) = 0.$ Also note that $X$ will
be in $B(0,R_{\nu(m+1)}(c))$ for all large enough $n,$ so there is
some $i \leq \ka_m$ such that $X \in B(x_i^m,r_{\nu(m)}(\ep)).$
So, if $m$ is large enough,
\[ r_{\nu(m)}(\ep) + r_{\nu(m+1)}(u) \leq r_{\nu(m+1)}(t) \leq r_n(t).\]
So, $F_m(x_i(m))$ and hence $G_m$ occur. Since $G_m$ occurs
finitely often a.s., $d_n \leq r_n(t)$ for all large $n,$ a.s.
The result now follows since $\ep > 0$ and $c > 2$ are arbitrary.

Now we derive a lower bound for $d_n$. Let $r_n(t) = \frac{t
\log_2 n}{\sqrt{log n}}.$

\begin{prop} Let $t < \frac{d}{\sqrt{2}}.$ Then with probability
$1$, $d_n \geq r_n(t),$ eventually. \lab{prop2} \end{prop}

{\bf Proof. } We prove the above proposition using the
Poissonization technique, which uses the following Lemma (see
Lemma 1.4, Penrose (2003)).

\begin{lem} Let $N(\lam)$ be Poisson random variables with mean
$\lam.$ Then there exists a constant $c_1$ such that for all
$\lam > \lam_1,$
\[P[X > \lam+\lam^{3/4}/2]\leq c_1\exp(-\lam^{1/2}), \]
and
\[P[X < \lam-\lam^{3/4}/2]\leq c_1\exp(-\lam^{1/2}).\] \lab{l3}
\end{lem}

Enlarging the probability space, assume that for each $n$ there
exist Poisson variables $N(n)$ and $M(n)$ with means $n-n^{3/4}$
and $2n^{3/4}$ respectively, independent of each other and of
$\{X_1,X_2,\ldots \}.$ Define the point processes
\[ \cP_n^- = \{X_1,X_2,\ldots,X_{N(n)} \}, \qquad \cP_n^+ = \{
X_1,X_2,\ldots , X_{N(n)+M(n)} \}. \]
Then, $\cP_n^-$ and $\cP_n^+$ are Poisson point processes on
$\cR^d$ with intensity functions $(n-n^{3/4})\phi(\cdot)$ and
$(n+n^{3/4})\phi(\cdot)$ respectively. The point processes
$\cP_n^-$, $\cP_n^+$ and $\cX_n$ are coupled in such a way that
$\cP_n^- \subset \cP_n^+$. Thus, if $H_n = \{\cP_n^- \subset \cX_n
\subset \cP_n^+ \},$ then by the Borel-Cantelli Lemma and
Lemma~\ref{l3}, $P[H_n^c \mbox{ i.o. }] = 0.$ Hence $\{\cP_n^-
\subset \cX_n \subset \cP_n^+ \}$ a.s. for all large enough $n.$

Pick numbers $(2d+c-2)/2 \sqrt{2} < t < u.$ Let $\ep > 0$ satisfy
$\ep + t < u.$

Consider the annulus $A_n = B(0,R_n(c) \setminus
B(0,R_n^{\prime}),$ where $R_n(c)$ is as defined in (\ref{e4a})
and $R_n^{\prime} = R_n(-2).$ For each $n,$ choose a non-random
set $\{ x_1^n,x_2^n,\ldots ,x_{\sg_n}^n \} \subset A_n,$ such
that the balls $B(x_i^n,r_n(u)),$ $1 \leq i \leq \sg_n$ are
disjoint. The packing number $\sg_n$ is the maximum number of
disjoint balls $B(x,r_n(u)$, with $x \in A_n.$
\eq{ \sigma_n  \stackrel{>}{\sim}
\frac{R_n^d(c)-{R^{\prime}_n}^d}{r_n^d(u)} \sim \left(\frac{\log
\; n}{\log_2 n}\right)^{d-1}. \lab{p2}}

Let $d_n^o$ be the LNND of the points of $\cX_n$ that fall in
$A_n$. By Lemma~\ref{l2}, there will be points in $A_n$ for all
large enough $n$, a.s. For any point process ${\cal{X}}$ and any
$B\subset \Re^d,$ let ${\cal{X}}[B]$ be the number of $X$ in $B.$
Let $E_n(x)$ be the event such that
\[E_n(x) = \{{\cal{P}}^{-}_n[B(x,r_n(\ep))]= 1\}
\cap \{{\cal{P}}^{+}_n[B(x,r_n(u))]= 1\}, \]
where $x\in B(0,R_n(c)) \setminus B(0,R_n^{\prime}).$
Set ${\cal{I}}_n = {\cal{P}}_{n}^{+}\setminus {\cal{P}}_{n}^{-},$
and set $U_n(x) = B(x,r_n(\ep)),$ and $V_n(x) =
B(x,r_n(u))\setminus U_n(x),$ then for each $n$ and $x$ the random
variables ${\cal{P}}_{n}^{-}(U_n), {\cal{P}}_{n}^{-}(V_n),
{\cal{I}}_n(V_n),$ and ${\cal{I}}_n(U_n),$ are independent
Poissons, and $E_n(x)$ is the event that the first of these
variables is 1 and the others are zero. Thus,
\begin{eqnarray}
P[E_n] & = & (n-n^{3/4})\int_{B(x,r_n(\ep))}\phi(y)dy
\exp\left(-(n+n^{3/4})\int_{B(x,r_n(u))}\phi(y)dy\right)\nonumber\\
& = & (n-n^{3/4})I(x,r_n(\ep))
\exp\left(-(n+n^{3/4})I(x,r_n(u))\right)\nonumber\\
& \geq & (n-n^{3/4})I(R_n^{\prime},r_n(\ep))
\exp\left(-(n+n^{3/4})I(R_n(c),r_n(u))\right)\nonumber\\
& \sim & C_d n r_n^d(\ep)
\exp(R_n^{\prime}r_n(\ep)-{R_n^{\prime}}^2/2)
(R_n^{\prime}r_n(\ep))^{-\frac{d+1}{2}}\nonumber\\
&& \cdot \exp\left(-C_d n r_n^d(u) \exp(R_n(c)r_n(u)-R^{2}_n(c)/2)
(R_n(c)r_n(u))^{-\frac{d+1}{2}}\right) \nonumber\\
& \sim & \frac{(\log _2\; n)^{\frac{d-1}{2}}} {(\log \;
n)^{d-2-\ep\sqrt{2}}}\exp\left(-\frac{C_d 2^{-\frac{d+1}{4}}
(u\log _2\; n)^{\frac{d-1}{2}}}{(\log \; n)^{d+c/2-1-u\sqrt{2}}}\right)\nonumber\\
& \sim & \frac{(\log _2\; n)^{\frac{d-1}{2}}} {(\log \;
n)^{d-2-\ep\sqrt{2}}}, \label{proba of en}
\end{eqnarray}
where the last relation follows since $u <
\frac{2d+c-2}{2\sqrt{2}}.$

If $H_n$ and $E_n(x)$ happen, then there is a point $X \in \cX_n$
in $B(0,R_n(c))\setminus B(0,R_n^{\prime})$ with no other point of
${\cal{X}}_n$ in $B(x, r_n(u)).$ Therefore
\begin{equation}
\{d_n^o \leq r_n\} \subset H_n^c \cup
\left(\bigcup_{i=1}^{\sigma_n}E_n(x_i^n)\right)^c.\label{dn leq
rn}
\end{equation}
From above, $P[H_n^c]$ is summable. We will show that
$P[\cup_{i=1}^{\sigma_n}E_n(x_i^n)]$ is summable. Since by
Lemma~\ref{l2}, there are points of $\cX_n$ in $A_n$ infinitely
often, a.s., we conclude that $d_n^o$ and hence $d_n$ will be
greater than $r_n(t)$ infinitely often a.s.

The events $E_n(x_i^n),\:\:1\leq i\leq \sigma_n$ are independent,
so by (\ref{proba of en}), for large enough $n$,
\begin{eqnarray*}
P\left[\left(\bigcup_{i=1}^{\sigma_n}E_n(x_i^n)\right)^c\right]
& \leq & \prod_{i=1}^{\sigma_n}\exp(-P[E_n(x_i^n)]) \\
& \leq & \exp\left(-C_1 \sigma_n\frac{(\log _2\;
n)^{\frac{d-1}{2}}}
{(\log \; n)^{d-2-\ep\sqrt{2}}}\right) \\
& \sim & \exp\left(-C_2 \left(\frac{\log \; n}{\log _2\;
n}\right)^{d-1} \frac{(\log _2\; n)^{\frac{d-1}{2}}}
{(\log \; n)^{d-2-\ep\sqrt{2}}}\right) \\
& = & \exp\left(- C_2 \frac{(\log \;
n)^{\ep\sqrt{2}+1}}{(\log_2{n})^{(d-1)/2}}\right),
\end{eqnarray*}
which is summable in $n$. $C_1$ and $C_2$ are some constants. This
completes the proof of Proposition~\ref{prop2}.

Theorem~\ref{t1} now follows from Propositions~\ref{prop1} and
\ref{prop2}.


\begin{thebibliography}{99}
\bibitem{Appel} Appel, M. J. B. and Russo, R.P. (1997), The minimum
vertex degree of a graph on the uniform points in $[0,1]^d$,
\textit{Advances in Applied Probability}, {\bf 29}, 582-594.
%
\bibitem{Henze} Dette, H. and Henze, N. (1989), The limit
distribution of the largest neighbor link in the unit d-cube,
\textit{Journal of Applied Probability}, {\bf 26}, 67-80.
%
\bibitem{Penrose1} Penrose, M. (1998), Extremes for the minimal
spanning tree on the Normally distributed points,
\textit{Advances in Applied Probability}, {\bf 30}, 628-639.
%
\bibitem{Penrose2} Penrose, M. (1999), A strong law for the largest
 nearest neighbor link between random points,
\textit{Journal of the london mathematical society}, {\bf 60},
951-960.
%
\bibitem{Penrose} Penrose, M. (2003), Random Geometric Graphs,
\textit{Oxford University Press.}
%
\bibitem{Steel} Steele, J. M. and Tierney, L. (1986), Boundary
dominaton and the distribution of the largest nearest-neighbor
link, \textit{Journal of Applied Probability}, {\bf 23}, 524-528.
\end{thebibliography}
\end{document}